\documentclass[12pt]{article}
\usepackage{amsmath}
\usepackage{graphicx,psfrag,epsf}
\usepackage{natbib}
\usepackage[margin=1in]{geometry}

\newtheorem{lemma}{Lemma}

\begin{document}

\bibliographystyle{natbib}

\title{\bf Anova of Balanced Variance Component Models}
  \author{Martin Bilodeau\hspace{.2cm}\\
    D\'epartement de math\'ematiques et de statistique\\ Universit\'e de Montr\'eal\\
martin.bilodeau@umontreal.ca}
\date{}
\maketitle

\bigskip
\begin{abstract}
Balanced linear models with fixed effects are taught in undergraduate programs of all universities. These occur in experimental designs such as one-way and two-way Anova, randomized complete block designs (RCBD) and split plot designs. The distributional theory for sums of squares in fixed effects models can be taught using the simplest form of Cochran's Theorem. The contribution provided here allows for an easy extension of the distributional theory to corresponding models with random effects. The main tool used is a simple result on noncentral chi-square distribution overlooked in textbooks at undergraduate and  graduate levels.  
\end{abstract}

\noindent%
\vfill

\newpage

\section{Introduction}
\label{sec:intro}

The one-way model with fixed effects is
\begin{equation}\label{eqn-onef}y_{ij}=\mu+\alpha_i+e_{ij},\ i=1,\dots,a,\  j=1,\dots,n,\end{equation}
where $\alpha_1$,\dots, $\alpha_a$ are fixed effects of the factor $A$ and $e_{ij}\sim N(0,\sigma^2)$ are error terms. 
Students are shown, see \textit{e.g.} \cite{MR2552961},
the Anova orthogonal decomposition of total sum of squares into a sum of squares due to factor $A$ and an error sum of squares,
\[SS_T=SS_A+SS_E.\]
The analysis usually proceeds by proving with Cochran's Theorem that these two sums of squares are independent and distributed as 
\begin{eqnarray*}
SS_A&\sim&\sigma^2\chi^2(a-1,\gamma_A/\sigma^2),\ \text{where } \gamma_A=n\sum_{i=1}^a(\alpha_i-\bar\alpha)^2\\
SS_E&\sim&\sigma^2\chi^2(a(n-1)),
\end{eqnarray*}
where $\gamma_A$ is called the noncentrality parameter of the noncentral chi-square distribution with $a-1$ degrees of freedom. 
The ratio then gives the $F$ test statistic for the hypothesis $H:\lambda_A=0$
\[F_A=\frac{MS_A}{MS_E}\sim F(a-1,a(n-1),\gamma_A/\sigma^2), \]
where $MS_A=SS_A/(a-1)$ and $MS_E=SS_E/[a(n-1)]$ are mean sums of squares.

Going from the one-way model with fixed effects to the model with random effects, the model becomes
\begin{equation}\label{eqn-oner}y_{ij}=\mu+\alpha_i+e_{ij},\ i=1,\dots,a,\  j=1,\dots,n,\end{equation}
where $\alpha_i\sim N(0,\sigma^2_\alpha)$ are random effects of the factor $A$ and $e_{ij}\sim N(0,\sigma^2)$ are errors. 
All random variables $\alpha_i$ and $e_{ij}$ are mutually independent.
The hypothesis related to factor $A$ is now $H:\sigma^2_\alpha=0$. 
A difficulty for the student to derive the distribution of $F_A$ is that observations sharing the same level of the factor $A$ are now correlated,
\begin{eqnarray*}
cov(y_{ij},y_{i'j'}) &=&0,\ i\neq i'\\
cov(y_{ij},y_{ij'}) &=&\sigma^2_\alpha,\ j\neq j'\\
var(y_{ij})&=&\sigma^2+\sigma^2_\alpha.
\end{eqnarray*}
The parameters $\sigma^2$ and $\sigma^2_\alpha$ are called components of variance. The one-way model with random effects is a variance component model.
For this model, the student in an undergraduate course is often served the laconic ``it can be shown" that the null distribution of $F_A$ remains the same, \textit{i.e.} 
$F_A\sim  F(a-1,a(n-1))$. At best, in textbooks at the undergraduate level,  expected mean squares $E(MS_A)$ and  $E(MS_E)$ are obtained as a justification.
This paper proposes a simple way to extend the complete description of the distributional results from the fixed effects model to the random effects situation.
The method requires only basic tools known to students from an undergraduate first course of probabilities: the moment generating function and conditional distributions.
The method is based on the following property of noncentral chi-square distributions. 

\begin{lemma}\label{lemma-random}If $x|\gamma_1\sim c_1\chi^2(p,\gamma_1/c_1)$ and $\gamma_1\sim c_2\chi^2(p,\gamma_2/c_2)$ then, 
\[x\sim(c_1+c_2)\chi^2(p,\gamma_2/(c_1+c_2)).\]
\end{lemma}

\noindent
Proof: Using the moment generating function,
\begin{eqnarray*}
E \exp(tx)&=&E E[\exp(tx)|\gamma_1]\\
&=&E E\left\{\exp\left[tc_1\chi^2(p,\gamma_1/c_1)\right]|\gamma_1\right\}\\
&=&(1-2tc_1)^{-p/2}E\exp\left[\gamma_1 t/(1-2tc_1)\right]\\
&=&(1-2tc_1)^{-p/2}E\exp\left[c_2\chi^2(p,\gamma_2/c_2) t/(1-2tc_1)\right]\\
&=&(1-2tc_1)^{-p/2}[1-2c_2 t/(1-2tc_1)]^{-p/2}\exp\left[\frac{\gamma_2 t/(1-2tc_1)}{1-2tc_2/(1-2tc_1)} \right]\\
&=&(1-2tc_1-2tc_2)^{-p/2}\exp\left[\gamma_2 t/(1-2tc_1-2tc_2) \right]\\
&=&[1-2t(c_1+c_2)]^{-p/2}\exp\left\{\left[\gamma_2/(c_1+c_2)\right] t(c_1+c_2)/\left[1-2t(c_1+c_2)\right]\right\}.
\end{eqnarray*}
In a first version of this lemma the variable $\gamma_1$ was assumed to be distributed as central chi-square. In a private communication, Eric Marchand informed the author of this more general result.
Lemma~\ref{lemma-random} can be applied to models listed in the abstract to easily derive the full distributional theory of all sums of squares adding to the total sum of squares.
\section{One-Way Anova}
\label{sec:one}

The joint distribution of $S=SS_A$ and $T=SS_E$ is now obtained in the one-way model with random effects. 
The first step is to recognize that, conditionally on the random effects $\alpha_i$,  model (\ref{eqn-oner}) reduces to model (\ref{eqn-onef}).
The conditional  joint density of $(S,T)$ given $ \gamma_A=n\sum_{i=1}^a(\alpha_i-\bar\alpha)^2$ is
\[f_{S,T}(s,t|\gamma_A)=f_{S}(s|\gamma_A) f_T(t).\]
Since $\sum_{i=1}^a(\alpha_i-\bar\alpha)^2\sim \sigma^2_\alpha\chi^2(a-1)$, it follows that $\gamma_A\sim n\sigma^2_\alpha\chi^2(a-1)$.
The unconditional joint density is obtained by integrating with respect to the density
$g(\gamma_A)$,
\[f_{S,T}(s,t)=\int_0^\infty f_{S}(s|\gamma_A)g(\gamma_A)d\gamma_A \cdot f_T(t).\]
The final result follows from Lemma~\ref{lemma-random} with
$p=a-1$, $c_1=\sigma^2$, $c_2=n\sigma^2_\alpha$ and $\gamma_2=0$.
The variables $SS_A$ and $SS_E$ are thus independent and distributed as
\begin{eqnarray*}
SS_A&\sim&(\sigma^2+n\sigma^2_\alpha)\chi^2(a-1)\\
SS_E&\sim&\sigma^2\chi^2(a(n-1)).
\end{eqnarray*}
The hypothesis can be tested using  
\[F_A=\frac{MS_A}{MS_E}\sim\frac{(\sigma^2+n\sigma^2_\alpha)}{\sigma^2}F(a-1,a(n-1)).\]

\section{Randomized Complete Block Design}
\label{sec:rcbd}

The RCBD model with fixed effects for the factor $A$ and the bloc factor $B$ is
\begin{equation}\label{eqn-rcbd}y_{ij}=\mu+\alpha_i+\beta_j+e_{ij},\ i=1,\dots,a,\ j=1,\dots,b,\end{equation}
where all the effects $\alpha_i$ and $\beta_j$ are fixed.
The Anova orthogonal decomposition is 
\[SS_T=SS_A+SS_B+SS_E\]
and, from Cochran's Theorem, the student is taught that
\begin{eqnarray*}
S=SS_A&\sim&\sigma^2\chi^2(a-1,\gamma_A/\sigma^2),\  \gamma_A=b\sum_{i=1}^a(\alpha_i-\bar\alpha)^2\\
T=SS_B&\sim&\sigma^2\chi^2(b-1,\gamma_B/\sigma^2),\  \gamma_B=a\sum_{j=1}^b(\beta_j-\bar\beta)^2\\
U=SS_E&\sim&\sigma^2\chi^2((a-1)(b-1) )
\end{eqnarray*}
are mutually independent. 

In the mixed model, the effects of factor $A$ are fixed and those of factor $B$ are random. 
The assumption is that all variables $\beta_j\sim N(0,\sigma^2_\beta)$ and
$e_{ij}\sim N(0,\sigma^2)$ are mutually independent. The hypothesis on factor $A$ is $H:\lambda_A=0$ and that for the bloc factor $B$ is $H:\sigma^2_\beta=0$.
Conditionally on the random effects $\beta_j$,  the RCBD model with random effects reduces to model (\ref{eqn-rcbd}) .
The conditional  joint density of $(S,T,U)$ given $\gamma_B=a\sum_{j=1}^b(\beta_j-\bar\beta)^2$  is
\[f_{S,T,U}(s,t,u|\gamma_B)=f_S(s)\cdot f_{T}(t|\gamma_B) f_U(u).\]
Since $\sum_{j=1}^b(\beta_j-\bar\beta)^2\sim \sigma^2_\beta \chi^2(b-1)$, it follows that $\gamma_B\sim a\sigma^2_\beta\chi^2(b-1)$.
The unconditional joint density is obtained by integrating with respect to the density
$g(\gamma_B)$,
\[f_{S,T,U}(s,t,u)=f_S(s)\cdot\int_0^\infty f_{T}(t|\gamma_B)g(\gamma_B)d\gamma_B\cdot f_U(u).\]
The final result follows from Lemma~\ref{lemma-random} with
$p=b-1$, $c_1=\sigma^2$, $c_2=a\sigma^2_\beta$ and $\gamma_2=0$.
The variables $SS_A$, $SS_B$ and $SS_E$ are thus mutually independent and distributed as
\begin{eqnarray*}
SS_A&\sim&\sigma^2\chi^2(a-1,\gamma_A/\sigma^2),\  \gamma_A=b\sum_{i=1}^a(\alpha_i-\bar\alpha)^2\\
SS_B&\sim&(\sigma^2+a\sigma^2_\beta)\chi^2(b-1)\\
SS_E&\sim&\sigma^2\chi^2((a-1)(b-1)).
\end{eqnarray*}
The hypotheses can thus be tested using  
\begin{eqnarray*}
F_A&=&\frac{MS_A}{MS_E}\sim F(a-1,(a-1)(b-1),\gamma_A/\sigma^2)\\
F_B&=&\frac{MS_B}{MS_E}\sim \frac{(\sigma^2+a\sigma^2_\beta)}{\sigma^2}F(b-1,(a-1)(b-1)).\\
\end{eqnarray*}
The null distribution of both $F$ tests is thus the same as in the pure fixed effects model (\ref{eqn-rcbd}). The error sum of squares is now written for later use in Comment~\ref{comment-inter} of Section~\ref{sec:comments}
\[SS_E=\sum_{i=1}^a\sum_{j=1}^b (y_{ij}-\bar y_{i.}-\bar y_{.j}+\bar y_{..})^2\sim\sigma^2\chi^2((a-1)(b-1)).\]

\section{Comments}
\label{sec:comments}

\begin{enumerate}

\item Lemma~\ref{lemma-random} on noncentral chi-square distributions.

This distributional property was not found in the statistical litterature including research papers and textbooks. The closest related result is in
\cite{MR663453} which considers improved estimation of the noncentrality parameter of chi-square distributions. They used the chi-square prior distribution and derived the posterior
mean $E(\lambda|x)$. The marginal distribution of $x$ was not given since it is not required to compute the Bayes estimator.

\item Cochran's Theorem. 

For  balanced linear models with only fixed effects, it was mentioned that Cochran's Theorem is an expedient way of proving that sums of squares adding to the total sum of squares, $SS_T$, are mutually independent noncentral chi-square variables. However, for balanced variance component models, Cochran's Theorem is not applicable since  $SS_T$ is no longer chi-square but a linear combination of mutually independent chi-square variables.

\item\label{comment-inter} Random effects of  interaction.

The two-way and split plot variance component models have random effects of interaction which can be treated similarly. The sum of squares for interaction between factors $A$ and $B$ in the two-way model with pure fixed effects
\begin{equation}y_{ijk}=\mu+\alpha_i+\beta_j+(\alpha\beta)_{ij}+e_{ijk},\ i=1,\dots,a,\ j=1,\dots,b,\ k=1,\dots,n,\end{equation}
is distributed as 
\[SS_{AB}\sim\sigma^2\chi^2((a-1)(b-1),\gamma_{AB}/\sigma^2),\]
where
\[\gamma_{AB}=n\sum_{i=1}^a\sum_{j=1}^b [(\alpha\beta)_{ij}-\overline{(\alpha\beta)}_{i.}-\overline{(\alpha\beta)}_{.j}+\overline{(\alpha\beta)}_{..}]^2.\]
In the model with random effects of interaction, $(\alpha\beta)_{ij}\sim N(0, \sigma^2_{\alpha\beta})$, this noncentrality parameter is distributed as the error sum of squares of an RCBD model. Hence, 
\[\gamma_{AB}\sim n \sigma^2_{\alpha\beta}\chi^2((a-1)(b-1)),\]
from which it follows from Lemma~\ref{lemma-random} that
\[SS_{AB}\sim(\sigma^2+n \sigma^2_{\alpha\beta})\chi^2((a-1)(b-1)).\]

\item Quadratic forms.

The theory of quadratic forms in normal variates can of course be used in the following way, see \textit{e.g.} \cite{MR1461543}. Denote the various sums of squares comprising $SS_T$ as 
$q_i=Y^\top A_iY$, $i=1,\dots,r$. When matrices $A_i$ are idempotent, quadratic forms can be written as $q_i=\|A_iY\|^2$, where $Y$ is some multivariate normal with mean $\mu$ and covariance matrix $V$. Hence quadratic forms $q_i$ are mutually independent chi-square if and only if $A_iV$ are idempotent and $A_iVA_j=0$, $i\neq j$. This method of proof does not build from the known results of the purely fixed effects case. It requires results of the theory of quadratic forms in normal variates generally not studied in a course at the undergraduate level. Moreover, the matrix $V$ is expressed using the Kronecker product which has generally not been studied in an undergraduate course of linear algebra.

\item Constraints.

In most textbooks, in order for the normal equations to have a unique solution, constraints are  imposed on  fixed effects models such as 
$\sum_{i=1}^a \alpha_i=0$ in the one-way model. In this case, the noncentrality parameter associated with $SS_A$ becomes $\gamma_A/\sigma^2=n\sum_{i=1}^a \alpha_i^2/\sigma^2$. The distribution of 
$\gamma_A$ then requires the conditional distribution of $(\alpha_1,\dots,\alpha_a)|\sum_{i=1}^a \alpha_i=0$ which is $N(0, \sigma^2_\alpha(I_a-\frac{J_a^\top J_a}{a}))$, where $J_a$ is a vector of ones of dimension $a$. Now, since $I_a-\frac{J_a^\top J_a}{a}$ is an idempotent matrix of rank $a-1$, $\sum_{i=1}^a \alpha_i^2|\sum_{i=1}^a \alpha_i=0\sim \sigma_\alpha^2\chi^2(a-1)$ as in Section~\ref{sec:one}. Such constraints were not used in this paper to avoid once again the theory of quadratic forms and the conditional distribution for multivariate normal distributions.

\end{enumerate}

\section{Conclusion}
\label{sec:conc}
Distributional results for balanced variance component models can be derived easily in an undergraduate course using Lemma~\ref{lemma-random} on noncentral chi-square distributions. This requires 
very little additional work once the results for the pure fixed effects model are obtained. This was demonstrated on two models: the one-way model and the RCBD model.
In exactly the same manner, the student can himself prove results for other mixed models such as the two-way and split plot models, once he or she has been taught distributional results for  pure effects models. All the examples required only the case $\gamma_2=0$ in  Lemma~\ref{lemma-random}. A detailed examination of the split-plot model shows that the extended case $\gamma_2>0$ is needed.


\begin{thebibliography}{}

\bibitem[\protect\citeauthoryear{Montgomery}{Montgomery}{2019}]{MR2552961}
Montgomery, D.~C. (2019).
\newblock {\em Design and analysis of experiments\/} (Tenth ed.).
\newblock John Wiley \& Sons, Inc., Hoboken, NJ.

\bibitem[\protect\citeauthoryear{Saxena and Alam}{Saxena and
  Alam}{1982}]{MR663453}
Saxena, K. M.~L. and K.~Alam (1982).
\newblock Estimation of the noncentrality parameter of a chi squared
  distribution.
\newblock {\em Ann. Statist.\/}~{\em 10\/}(3), 1012--1016.

\bibitem[\protect\citeauthoryear{Searle}{Searle}{1997}]{MR1461543}
Searle, S.~R. (1997).
\newblock {\em Linear models}.
\newblock Wiley Classics Library. John Wiley \& Sons, Inc., New York.
\newblock Reprint of the 1971 original.

\end{thebibliography}
\end{document}